\documentclass[a4paper]{amsart}

\usepackage{amsfonts}

\title[Simultaneous avoidance of generalized patterns]{Simultaneous avoidance of generalized patterns}
\author{Sergey Kitaev and Toufik Mansour}
\address{{\small Matematik, Chalmers tekniska h\"ogskola och G\"oteborgs
universitet, 412~96  G\"oteborg, Sweden}}
\email{kitaev@math.chalmers.se} 
\address{{\small LaBRI, Universit\'e Bordeaux 1, 351 cours de la Lib\'eration, 33405 Talence Cedex, France}}
\email{toufik@labri.fr} 

\newtheorem{prop}{Proposition}
\newtheorem{lemma}[prop]{Lemma}

\newtheorem{thm}[prop]{Theorem}
\theoremstyle{definition}

\newtheorem{remark}[prop]{Remark}

\begin{document}
\begin{abstract}
In \cite{BabStein} Babson and Steingr\'{\i}msson introduced
generalized permutation patterns that allow the requirement
that two adjacent letters in a pattern must be adjacent in
the permutation. In~\cite{Kit1} Kitaev considered simultaneous
avoidance (multi-avoidance) of two or more 3-patterns with no
internal dashes, that is, where the patterns correspond to
contiguous subwords in a permutation. There either an
explicit or a recursive formula was given for all but
one case of simultaneous avoidance of more than two patterns.

In this paper we find the exponential generating function for the
remaining case. Also we consider permutations that avoid a pattern
of the form $x-yz$ or $xy-z$ and begin with one of the patterns
$12\ldots k$, $k(k-1)\ldots 1$, $23\ldots k1$, $(k-1)(k-2)\ldots
1k$ or end with one of the patterns $12\ldots k$, $k(k-1)\ldots
1$, $1k(k-1)\ldots 2$, $k12\ldots (k-1)$. For each of these cases
we find either the ordinary or exponential generating functions or
a precise formula for the number of such permutations. Besides we
generalize some of the obtained results as well as some of the
results given in~\cite{Kit3}: we consider permutations avoiding
certain generalized 3-patterns and beginning (ending) with an
arbitrary pattern having either the greatest or the least letter
as its rightmost (leftmost) letter.
\end{abstract}

\maketitle

\thispagestyle{empty}

\section{Introduction and Background}
{\bf Permutation patterns:} All permutations in this
paper are written as words $\pi=a_1 a_2\ldots a_n$, where the
$a_i$ consist of all the integers $1,2,\ldots,n$. Let $\alpha\in S_n$ and $\tau\in S_k$
be two permutations. We say that $\alpha$ {\em contains} $\tau$ if
there exists a subsequence $1\leq i_1<i_2<\cdots<i_k\leq n$ such
that $(\alpha_{i_1},\dots,\alpha_{i_k})$ is order-isomorphic to
$\tau$, that is, for all $j$ and $m$, ${\tau}_j < {\tau}_m$ if and only if $a_{i_j} < a_{i_m}$; in such a context $\tau$ is usually called a {\em
pattern}. We say that $\alpha$ {\em avoids} $\tau$, or is $\tau$-{\em
avoiding}, if $\alpha$ does not contain $\tau$. The set of all
$\tau$-avoiding permutations in $S_n$ is denoted by $S_n(\tau)$. For
an arbitrary finite collection of patterns $T$, we say that
$\alpha$ avoids $T$ if $\alpha$ avoids each $\tau\in T$; the
corresponding subset of $S_n$ is denoted by $S_n(T)$.

While the case of permutations avoiding a single pattern has
attracted much attention, the case of multiple pattern avoidance
remains less investigated. In particular, it is natural, as the
next step, to consider permutations avoiding pairs of patterns
$\tau_1$, $\tau_2$. This problem was solved completely for
$\tau_1,\tau_2\in S_3$ (see \cite{SchSim}), for $\tau_1\in S_3$
and $\tau_2\in S_4$ (see \cite{W}), and for $\tau_1,\tau_2\in S_4$
(see \cite{B1,Km} and references therein). Several recent papers
\cite{CW,MV1,Kr,MV3,MV2} deal with the case $\tau_1\in S_3$,
$\tau_2\in S_k$ for various pairs $\tau_1,\tau_2$.

\medskip

{\bf Generalized permutation patterns:} In \cite{BabStein} Babson
and Steingr\'{\i}msson introduced {\em generalized permutation
patterns (GPs)} where two adjacent letters in a pattern may be
required to be adjacent in the permutation. Such an adjacency
requirement is indicated by the absence of a dash between the
corresponding letters in the pattern. For example, the permutation
$\pi=516423$ has only one occurrence of the pattern $2$-$31$,
namely the subword 564, but the pattern $2$-$3$-$1$ occurs also in
the subwords 562 and 563. Note that a classical pattern should, in
our notation, have dashes at the beginning and end. Since most of
the patterns considered in this paper satisfy this, we suppress
these dashes from the notation. Thus, a pattern with no dashes
corresponds to a contiguous subword anywhere in a permutation.
 The motivation for introducing these patterns was the study
of Mahonian statistics. A number of results on GPs were obtained
by Claesson, Kitaev and Mansour. See for example \cite{Claes},
\cite{Kit1, Kit2, Kit3} and \cite{Mans1, Mans2, Mans3}.

As in \cite{SchSim}, dealing with the classical patterns, one can
consider the case when permutations have to avoid two or more
generalized patterns simultaneously. A complete solution for
the number of permutations avoiding a pair of 3-patterns of
type (1,2) or (2,1), that is, the patterns having one internal
dash, is given in~\cite{ClaesMans}. In~\cite{Kit1} Kitaev gives
either an explicit or a recursive formula for all but
one case of simultaneous avoidance of more than two patterns.
 This is the case of avoiding the GPs 123, 231 and 312
simultaneously. In Theorem~\ref{avoidance_123_231_312}
we find the exponential generating function (e.g.f.)
for the number of such permutations.

As it was discussed in~\cite{Kit3}, if a permutation begins (resp.
ends) with the pattern $p = p_1p_2\ldots p_k$, that is, the $k$
leftmost (resp. rightmost) letters of the permutation form the
pattern $p$, then this is the same as avoidance of $k!-1$ patterns
simultaneously. For example, beginning with the pattern $123$ is
equivalent to the simultaneous avoidance of the patterns $(132]$,
$(213]$, $(231]$, $(312]$ and $(321]$ in the
Babson-Steingr\'{\i}msson notation. Thus demanding that a
permutation must begin or end with some pattern, in fact, we are
talking about simultaneous avoidance of generalized patterns. The
motivation for considering additional restrictions such as
beginning or ending with some patterns is their connection to some
classes of trees. An example of such a connection can be found
in~\cite[Theorem 5]{Kit3}. There it was shown that there is a
bijection between $n$-permutations avoiding the pattern 132 and
beginning with the pattern 12 and {\em increasing rooted trimmed
trees} with $n+1$ nodes. We recall that a trimmed tree is a tree
where no node has a single leaf as a child (every leaf has a
sibling) and in an increasing rooted tree, nodes are numbered and
the numbers increase as we move away from the root. The avoidance
of a generalized 3-pattern $p$ with no dashes and, at the same
time, beginning or ending with an increasing or decreasing pattern
was discussed in~\cite{Kit3}. Theorem~\ref{generalization_1}
generalizes some of these results to the case of beginning (resp.
ending) with an arbitrary pattern avoiding $p$ and having the
greatest or least letter as the rightmost (resp. leftmost) letter.

Propositions 4 -- 15 (resp. 16 -- 27) give a complete description
for the number of permutations avoiding a pattern of the form
$x-yz$ or $xy-z$ and beginning with one of the patterns $12\ldots
k$ or $k(k-1)\ldots 1$ (resp. $23\ldots k1$ or $(k-1)(k-2)\ldots
1k$). For each of these cases we find either the ordinary or
exponential generating functions or a precise formula for the
number of such permutations. Theorem~\ref{generalization_2}
generalizes some of these results. Besides, the results from
Propositions 4--27 give a complete description for the number of
permutations that avoid a pattern of the form $x-yz$ or $xy-z$ and
end with one of the patterns $12\ldots k$, $k(k-1)\ldots 1$,
$1k(k-1)\ldots 2$ and $k12\ldots (k-1)$. To get the last one of
these we only need to apply the reverse operation discussed in the
next section. The results of Theorems~\ref{generalization_1}
and~\ref{generalization_2} can also be used to get the case of
ending with a pattern from the sets $\Delta_k^{min}$ or
$\Delta_k^{max}$ introduced in the next section.

Except for the empty permutation, every permutation ends and
begins with the pattern $p=1$. To simplify the discussion we
assume that the empty permutation also begin with the pattern $1$.
This does not course any harm since, to count the generating
functions in question for this, we need only subtract $1$ from the
generating functions obtained in this paper.

\section{Preliminaries}\label{Preliminaries}
The {\em reverse} $R(\pi)$ of a permutation $\pi=a_1a_2 \ldots
a_n$ is the permutation $a_n \ldots a_2a_1$. The {\em complement}
$C(\pi)$ is the permutation $b_1b_2 \ldots b_n$ where
$b_i=n+1-a_i$. Also, $R \circ C$ is the composition of $R$ and
$C$. For example, $R(13254)=45231$, $C(13254)=53412$ and $R \circ
C(13254)=21435$. We call these bijections of $S_n$ to itself {\em
trivial}, and it is easy to see that for any pattern $p$ the
number $A_p(n)$ of permutations avoiding the pattern $p$ is the
same as for the patterns $R(p)$, $C(p)$ and $R \circ C(p)$. For
example, the number of permutations that avoid the pattern 132 is
the same as the number of permutations that avoid the pattern 231.
This property holds for sets of patterns as well. If we apply one
of the trivial bijections to all patterns of a set $G$, then we
get a set $G^{\prime}$ for which $A_{G^{\prime}}(n)$ is equal to
$A_G(n)$. For example, the number of permutations avoiding $\{
123, 132 \}$ equals the number of those avoiding  $\{ 321, 312\}$
because the second set is obtained from the first one by
complementing each pattern.

In this paper we denote the $n$th Catalan number by $C_n$; the
generating function for these numbers by $C(x)$; the $n$th Bell
number by $B_n$.

Also, $N^{p}_{q}(n)$ denotes the number of permutations that avoid
the pattern $q$ and begin with the pattern $p$; $G^{p}_{q}(x)$
(resp. $E^{p}_{q}(x)$) denotes the ordinary (resp. exponential)
generating function for the number of such permutations. Besides,
$\Gamma_k^{min}$ (resp. $\Gamma_k^{max}$) denotes the set of all
$k$-patterns with no dashes such that the least (resp. greatest)
letter of a pattern is the rightmost letter; $\Delta_k^{min}$
(resp. $\Delta_k^{max}$) denotes the set of all $k$-patterns with
no dashes such that the least (resp. greatest) letter of a pattern
is the leftmost letter.

Recall the following properties of $C(x)$:
\begin{equation}
C(x) = \frac{1-\sqrt{1-4x}}{2x} = \frac{1}{1-xC(x)}.
\label{catalan}
\end{equation}
\section{Simultaneous avoidance of $123$, $231$ and $312$}
The {\em Entringer numbers} $E(n,k)$ (see \cite[Sequence A000111/M1492]{SloPlo})
are the number of permutations on $1,2,\dots,n+1$,
starting with $k+1$, which, after initially falling, alternately
fall then rise. The Entringer numbers (see \cite{Ent}) are given
by
    $$E(0,0)=1,\qquad E(n,0)=0,$$
together with the recurrence relation
                $$E(n,k)=E(n,k+1)+E(n-1,n-k).$$
The numbers $E(n)=E(n,n)$,
are the secant and tangent numbers given by the generating
function
            $$\sec x+\tan x.$$

The following theorem completes the consideration of
multi-avoidance of more than two generalized 3-patterns with no
dashes made in~\cite{Kit1}.

\begin{thm}\label{avoidance_123_231_312} Let $E(x)$ be the e.g.f.
for the number of permutations that avoid $123$, $231$ and $312$ simultaneously. Then
        $$E(x) = 1 + x(\sec(x) + \tan(x)).$$
\end{thm}
\begin{proof}
Let $s(n;i_1,\dots,i_m)$ denote the number of permutations $\pi\in S_n(123,231,312)$
such that $\pi_1\pi_2\dots\pi_m=i_1i_2\dots i_m$ and $f: S_n\rightarrow S_n$ be a map defined by
$$f(\pi_1\pi_2\dots\pi_n)=(\pi_1+1)(\pi_2+1)\dots(\pi_n+1),$$
where the addition is modulo $n$. Using $f$ one can see that for all $a=1,2,\dots,n-1$,
\begin{equation}
s(n;a)=s(n;a+1).
\label{equation0}
\end{equation}
Thus, $|S_n(123,231,312)|=n s(n;1)$ and we only need to prove that $s(n;1)=E_{n-1}$, where $E_n$ is the $n$th Euler number (see \cite[Sequence A000111/M1492]{SloPlo}).

Suppose $\pi\in S_n(123,231,312)$ is an $n$-permutation such that $\pi_1=1$ and $\pi_2=t$.
Since $\pi$ avoids $123$, we get $\pi_3\leq t-1$ and it is easy to see that
$$s(n;1,t)=\sum_{j=2}^{t-1}s(n;1,t,j)=\sum_{j=1}^{t-2}s(n-1;t-1,j),$$
so
$$s(n;1,t+1)=s(n;1,t)+\sum_{j=1}^{t-1}s(n-1;t,j)-\sum_{j=1}^{t-2}s(n-1;t-1,j).$$
Using (\ref{equation0}) we get
$$s(n;1,t+1)=s(n;1,t)+s(n-1;t,1)+\sum_{j=2}^{t-1}s(n-1;t-1,j-1)-\sum_{j=1}^{t-2}s(n-1;t-1,j),$$
and by (\ref{equation0}) again, we have for all $t=2,3,\dots,n-1$,
        $$s(n;1,t+1)=s(n;1,t)+s(n-1;1,n-t+1).$$
Besides, by the definition, it is easy to see that $s(n;1,2)=0$ for all
$n\geq 3$, hence using the definition of Entringer numbers \cite{Ent} we get
$s(n;1)=\displaystyle\sum_{t=2}^n s_{n;1,t}=E_{n-1}$, as required.
\end{proof}

\section{Avoiding a 3-pattern with no dashes and beginning with a pattern whose rightmost letter is the greatest or smallest}

The following theorem generalizes Theorems 7 and 8 in~\cite{Kit3}. Recall the definition of $E_{q}^{p}$ in Section~\ref{Preliminaries}.

\begin{thm}\label{generalization_1} Suppose $p_1, p_2 \in \Gamma_k^{min}$ and $p_1 \in S_k(132)$, $p_2 \in S_k(123)$. Thus, the complements $C(p_1), C(p_2) \in \Gamma_k^{max}$ and $C(p_1) \in S_k(312)$, $C(p_2) \in S_k(321)$. Then, for $k \ge 2$,
$$E^{p_1}_{132}(x) = E^{C(p_1)}_{312}(x) = \frac{\int_{0}^{x}t^{k-1}e^{-t^2/2}\ dt}{(k-1)!(1 - \int_{0}^{x}e^{-t^2/2}\ dt)}$$
and
$$E^{p_2}_{123}(x) = E^{C(p_2)}_{321}(x) = \frac{e^{x/2}\int_{0}^{x}e^{-t/2}t^{k-1}\sin(\frac{\sqrt{3}}{2}t + \frac{\pi}{6}))\ dt}{(k-1)!\cos(\frac{\sqrt{3}}{2}x + \frac{\pi}{6})}.$$
\end{thm}

\begin{proof} To prove the theorem, it is enough to copy the proofs of  Theorems 7
and 8 in~\cite{Kit3}, since the fact that the first $k-1$ letters of $p$ are possibly
not in decreasing order is immaterial for the proofs of that theorems. Thus we can get the
formula for $E^{p}_{132}(x)$ and  $E^{p}_{123}(x)$, and automatically, using properties of the complement, the formula for $E^{C(p)}_{312}(x)$ and  $E^{C(p)}_{321}(x)$, directly from these theorems. However we give
here a proof of the formula for $E^{p}_{132}(x)$ and refer to~\cite[Theorem 8]{Kit3} for a proof of the formula for $E^{p}_{123}(x)$.

If $k=1$, we have no additional restrictions, that is, we are dealing only with the avoidance of $132$ and,
according to \cite[Theorem 4.1]{ElizNoy} or \cite[Theorem 12]{Kit2},
$$E^{1}_{132}(x) = \frac{1}{1 - \int_{0}^{x}e^{-t^2/2}\ dt}.$$

Also, according to \cite[Theorem 6]{Kit3},
$$E^{12}_{132}(x) = \frac{e^{-x^2/2}}{1 - \int_{0}^{x}e^{-t^2/2}\ dt} - x - 1.$$

Let $R_{n, k}$ (resp. $F_{n,k}$) denote the number of $n$--permutations that
avoid the pattern $132$ and begin with a decreasing (resp. increasing) subword of
length $k > 1$ and let $\pi$ be such a permutation of length $n+1$. Suppose $\pi = \sigma 1 \tau$. If $\tau = \epsilon$ then, obviously, there are $R_{n,k}$
ways to choose $\sigma$. If $|\tau| = 1$, that is, 1 is in the second position from the right, then there are $n$ ways to choose the rightmost letter in $\pi$ and
we multiply this by $R_{k,n-1}$, which is the number of ways to choose $\sigma$.
If $|\tau| > 1$ then $\tau$ must begin with the pattern $12$, otherwise the
letter 1 and the two leftmost letters of $\tau$ form the pattern $132$, which
is forbidden. So, in this case there are
$\displaystyle\sum_{i \geq 0}{n \choose i}R_{i,k}F_{n-i,2}$ such permutations
with the right properties, where $i$ indicates the length of $\sigma$.
In the last formula, of course, $R_{i,k} = 0$ if $i<k$. Finally we have to consider the situation when 1 is in the $k$-th position. In this case we can choose the
letters of $\sigma$ in ${n \choose k-1}$ ways, write them in decreasing order
and then choose $\tau$ in $F_{n-k+1,2}$ ways. Thus
\begin{equation}
R_{n+1,k} = R_{n,k} + nR_{n-1,k} +
\displaystyle\sum_{i\geq 0}{n \choose i}R_{i,k}F_{n-i,2} + {n \choose k-1}F_{n-k+1,2}.
\label{equality3}
\end{equation}

We observe that (\ref{equality3}) is not valid for $n=k-1$ and $n=k$.
Indeed, if 1 is in the $k$-th position in these cases, the term
${n \choose k-1}F_{n-k+1,2}$, which counts the number of such permutations, is zero, whereas there is one ``good''
$(n+1)$--permutation in case $n=k-1$ and $n$ good
$(n+1)$--permutations in case $n=k$. Multiplying both sides of the
equality with $x^n/n!$, summing over $n$ and
using the observation above (which gives the term $x^{k-1}/(k-1)! + kx^k/k!$
in the right-hand side of equality~(\ref{equality4})), we get
\begin{equation}
\frac{d}{dx}E^{p}_{132}(x) = (E^{12}_{132}(x) + x + 1)E^{p}_{132}(x) + (E^{12}_{132}(x)+x+1)\frac{x^{k-1}}{(k-1)!},
\label{equality4}
\end{equation}
with the initial condition $E^{p}_{132}(0) = 0$. We solve this equation and get
$$E^{p}_{132}(x) = \frac{E^{1}_{132}(x)}{(k-1)!}\int_{0}^{x}\frac{(E^{12}_{132}(t) + t + 1)t^{k-1}}{E^{1}_{132}(t)}\ dt = \frac{E^{1}_{132}(x)}{(k-1)!}\int_{0}^{x}t^{k-1}e^{-t^2/2}\ dt.$$
\end{proof}

\begin{remark} It is obvious that if in the previous theorem $p_1 \not\in S_k(132)$ and $p_2 \not\in S_k(123)$, then $E^{p_1}_{132}(x) = E^{p_2}_{123}(x) = 0$. \end{remark}
\section{Avoiding a pattern x-yz and beginning with an increasing or decreasing pattern}

In this section we consider avoidance of one of the patterns $1-23$, $1-32$, $2-31$, $2-13$, $3-12$ and $1-32$ and beginning with a decreasing pattern. We get all the other cases, that is, avoidance of one of these patterns and beginning with an increasing pattern, by the complement operation. For instance, we have $E^{k(k-1)\ldots 1}_{1-23}(x) = E^{12\ldots k}_{3-21}(x)$.

\begin{prop}\label{about_1-23} We have
\[ E_{1-23}^{k(k-1)\ldots 1}(x) = E_{1-32}^{k(k-1)\ldots 1}(x) = \left\{ \begin{array}{ll} (e^{e^x}/(k-1)!)\int_{0}^{x}t^{k-1}e^{-e^t+t}\ dt, & \mbox{if $k\geq 2$,}\\[2mm] e^{e^x-1}, & \mbox{if $k=1$.}
\end{array}
\right. \]
\end{prop}

\begin{proof} We prove the statement for the pattern $1-23$. All the arguments we give for this pattern are valid for the pattern $1-32$. The only difference is that instead of decreasing order in $\tau$ (see below), we have increasing order.

Suppose $k\geq 2$. Let $B_{n, k}$ denote the number of $n$--permutations that
avoid the pattern $1-23$ and begin with a decreasing subword of
length $k$. Suppose $\pi= \sigma 1 \tau$ be one of such permutations of length $n+1$. Obviously, the letters of $\tau$ must be in decreasing order since otherwise we have an occurrence of $1-23$ in $\pi$ starting from the letter 1. If $|\sigma| = i$ then we can choose the letters of $\sigma$ in ${n \choose i}$ ways. Since the letters of $\tau$ are in decreasing order, they do not affect $\sigma$ and thus there are $B_{i, k}$ possibilities to choose $\sigma$. Besides, if $|\sigma| = k-1$ and letters of $\sigma$ are in decreasing order, we get ${n \choose k-1}$ additional possibilities to choose $\pi$. Thus
$$B_{n+1,k} = \displaystyle\sum_{i\geq 0}{n \choose i}B_{i,k} + {n \choose k-1}.$$

Multiplying both sides of the
equality with $x^n/n!$ and summing over $n$, we get the differential equation
$$\frac{d}{dx}E_{1-23}^{k(k-1)\ldots 1}(x) = (E_{1-23}^{k(k-1)\ldots 1}(x) + \frac{x^{k-1}}{(k-1)!})e^x$$
with the initial condition $E_{1-23}^{k(k-1)\ldots 1}(0)=0$. The solution to this equation is given by
\begin{equation}
E_{1-23}^{k(k-1)\ldots 1}(x)=(e^{e^x}/(k-1)!)\int_{0}^{x}t^{k-1}e^{-e^t+t}\ dt.
\label{equality5}
\end{equation}

If $k=1$, then there is no additional restriction. According to \cite[Prop. 2]{Claes} (resp. \cite[Prop. 5]{Claes}), the number of $n$-permutations that avoid the pattern 1-23 (resp. 1-32) is the $n$th Bell number and the e.g.f. for the Bell numbers is $e^{e^x-1}$. However, all the arguments used for $k\geq 2$ remain the same for the case $k=1$ except for the fact that we do not count the empty permutation, which, of course, avoids 1-23. So, if $k=1$, we need to add 1 to the right-hand side of~(\ref{equality5}):
$$E_{1-23}^{1}(x)= e^{e^x}\int_{0}^{x}e^{-e^t+t}\ dt + 1 = e^{e^x-1}.$$
\end{proof}

\begin{prop}\label{about_3-12} We have
\[ E_{3-12}^{k(k-1)\ldots 1}(x) = \left\{ \begin{array}{ll} e^{e^x}\int_{0}^{x}e^{-e^t}\displaystyle\sum_{n\geq k-1}\frac{t^n}{n!}\ dt, & \mbox{if $k\geq 2$,}\\[2mm] e^{e^x-1}, & \mbox{if $k=1$.}
\end{array}
\right. \]
\end{prop}

\begin{proof} Suppose $k\geq 2$. Let $B_{n, k}$ denote the number of $n$--permutations that avoid the pattern $3-12$ and begin with a decreasing subword of
length $k$. Suppose $\pi = \sigma (n+1) \tau$ be such a permutation of length $n+1$. Obviously, the letters of $\tau$ must be in decreasing order since otherwise we have an occurrence of the pattern $3-12$ in $\pi$ starting from the letter $(n+1)$. If $|\sigma| = i$ then we can choose the letters of $\sigma$ in ${n \choose i}$ ways. Since the letters of $\tau$ are in decreasing order, they do not affect $\sigma$ and thus there are $B_{i, k}$ possibilities to choose $\sigma$. Besides, if $n\geq k-1$, then $\pi$ can be decreasing, that is, $(n+1)$ can be in the leftmost position. Thus
$$B_{n+1,k} = \displaystyle\sum_{i\geq 0}{n \choose i}B_{i,k} + \delta_{n,k},$$
where
\[ \delta_{n,k} = \left\{ \begin{array}{ll} 1, & \mbox{if $n\geq k-1$,}
\\ 0, & \mbox{else.}
\end{array}
\right. \]

Multiplying both sides of the equality with $x^n/n!$ and summing over $n$, we get the differential equation
$$\frac{d}{dx}E_{3-12}^{k(k-1)\ldots 1}(x) = e^xE_{3-12}^{k(k-1)\ldots 1}(x) + \displaystyle\sum_{n \geq k-1}\frac{x^n}{n!}$$
with the initial condition $E_{3-12}^{k(k-1)\ldots 1}(0)=0$. The solution to this equation is given by
\begin{equation}
E_{3-12}^{k(k-1)\ldots 1}(x) = e^{e^x}\int_{0}^{x}e^{-e^t}\displaystyle\sum_{n \geq k-1}\frac{t^n}{n!}\ dt.
\label{equality6}
\end{equation}

If $k=1$, then there is no additional restriction. In \cite[Prop. 5]{Claes} it is shown that $E_{1-32}^{1}(x)=e^{e^x-1}$. Using the complement, the number of $n$-permutations that avoid $1-32$ is equal to the number of $n$-permutations that avoid $3-12$. We get that $E_{3-12}^{1}(x)=e^{e^x-1}$. However, all the arguments used for the case $k\geq 2$ remain the same for the case $k=1$ except the fact that we do not count the empty permutation, which avoids 3-12. So, if $k=1$, we need to add 1 to the right-hand side of~(\ref{equality6}):
$$E_{3-12}^{1}(x)= e^{e^x}\int_{0}^{x}e^{-e^t}e^t\ dt + 1 = e^{e^x-1}.$$ \end{proof}

\begin{prop}\label{about_3-21} We have
\[ E_{3-21}^{k(k-1)\ldots 1}(x) = \left\{ \begin{array}{ll} 0, & \mbox{if $k\geq 3$,} \\[2mm] e^{e^x}\int_{0}^{x}e^{-e^t}(e^t - 1)\ dt, & \mbox{if $k=2$,}
\\[2mm] e^{e^x-1}, & \mbox{if $k=1$.}
\end{array}
\right. \]
\end{prop}

\begin{proof} For $k\geq 3$, the statement is obviously true. If $k=1$, then the statement follows from \cite[Prop. 2]{Claes} and the fact that there are as many $n$-permutations avoiding the pattern $1-23$, as $n$-permutations avoiding the pattern $3-21$. For the case $k=2$, we can use exactly the same arguments as those in the proof of Proposition~\ref{about_3-12} to get the same recurrence relation and thus the same formula, which, however, is valid only for $k=2$. \end{proof}

Recall the definition of $N_{q}^{p}$ in Section~\ref{Preliminaries}.

\begin{prop}\label{about_2-13} We have

\[ N_{2-13}^{k(k-1)\ldots 1}(n) = \left\{ \begin{array}{ll} C_{n-k+1}, & \mbox{if $n\geq k$,} \\ 0, & \mbox{else.}
\end{array}
\right. \]
\end{prop}

\begin{proof} If $k=1$, then the statement follows from~\cite[Prop. 22]{Claes}. Suppose $k\geq 2$ and let $\pi = \sigma n \tau$ be an $n$-permutation avoiding $2-31$ and beginning with the pattern $k(k-1)\ldots 1$. Suppose, without loss of generality that $\sigma$ consists of the letters $1, 2,\ldots, \ell$. Now $\ell$ must be the rightmost letter of $\sigma$, since otherwise $\ell$, the rightmost letter of $\sigma$ and $n$ form the pattern $2-13$. Also, the letter $(\ell - 1)$ must be next to the rightmost letter of $\sigma$ since otherwise the letter $(\ell - 1)$, next to the rightmost letter of $\sigma$ and the letter $\ell$ form the pattern $2-13$. And so on. Thus $\sigma$ must be increasing, which contradicts the fact that $\pi$ must begin with a decreasing pattern of length greater than 1. So $|\sigma| = 0$ and $\tau$ must begin with the pattern $(k-1)(k-2)\ldots 1$. Now, we can consider the letter $(n-1)$ and, by the same reasoning, get that it must be in the second position of $\pi$. Then we consider $(n-2)$, and so on up to the letter $(n-k+2)$. Finally, we get that $\pi = n(n-1)\ldots (n-k+2){\pi}^{\prime}$, where ${\pi}^{\prime}$ must avoid the pattern $2-13$ and thus, there are $C_{n-k+1}$ ways to choose $\pi$ (\cite[Prop. 22]{Claes}). \end{proof}

Recall that $C(x)$ is the generating function for the Catalan numbers. Also recall the definition of $G_{q}^{p}$ in Section~\ref{Preliminaries}.

\begin{prop}\label{about_2-31} We have
$$G^{k(k-1)\ldots 1}_{2-31}(x) = \left\{ \begin{array}{ll} x^kC^{k+1}(x), & \mbox{if $k\geq 2$} \\ C(x), & \mbox{if $k=1$.}\end{array}
\right.$$
\end{prop}

\begin{proof} If $k=1$, then there is no additional restriction, and thus $G^{1}_{2-31}(x) = C(x)$ (applying the complement operation to \cite[Prop. 22]{Claes}).

Suppose $k\geq  2$. Using the reverse, we see that beginning with $k(k-1)\ldots 1$ and avoiding $2-31$ is equivalent to ending with $12\ldots k$ and avoiding $13-2$, which by \cite{Claes} is equivalent to ending with $12\ldots k$ and avoiding $1-3-2$.

Suppose $\pi={\pi}^{\prime}n{\pi}^{\prime\prime}$ ends with $12\ldots k$ and avoids $1-3-2$. Each letter of ${\pi}^{\prime}$ must be greater than any letter of ${\pi}^{\prime\prime}$, since otherwise we have an occurrence of the pattern $1-3-2$ involving the letter $n$. Also, ${\pi}^{\prime}$ and ${\pi}^{\prime\prime}$ avoid the pattern $1-3-2$, and ${\pi}^{\prime\prime}$ ends with the pattern $12\ldots k$. In terms of generating functions (the generating function for the number of permutations ending with $12\ldots k$ and avoiding $1-3-2$ is, of course, $G^{k(k-1)\ldots 1}_{2-31}(x)$) this means that
\begin{equation}
G^{k(k-1)\ldots 1}_{2-31}(x) = xC(x)G^{k(k-1)\ldots 1}_{2-31}(x) + xG^{(k-1)\ldots 1}_{2-31}(x),
\label{equality15}
\end{equation}
where the rightmost term corresponds to the case when ${\pi}^{\prime\prime}$ is empty. Now, (\ref{catalan}) and~(\ref{equality15}) give
$$G^{k(k-1)\ldots 1}_{2-31}(x)=x^kC(x)/(1-xC(x))^k=x^kC^{k+1}(x).$$
\end{proof}

\section{Avoiding a pattern xy-z and beginning with an increasing or decreasing pattern}

First of all we state the following well-known binomial identity

\begin{equation}
\displaystyle\sum_{i=1}^{n-m-k+1}\binom{n-m-i}{k-1}\binom{m+i-1}{m}=\binom{n}{m+k}.
\label{equality8}
\end{equation}

Let $s_q(n)$ denote the cardinality of the set $S_n(q)$ and $s_q(n;i_1,i_2,\dots i_m)$ denote the number of permutations $\pi\in S_n(q)$ with $\pi_1\pi_2\dots\pi_m=i_1i_2\dots i_m$.

In this section we consider avoidance of one of the patterns 12-3, 13-2 and 23-1 and beginning with an increasing or decreasing pattern. We get all the other cases, which are avoidance of one of the patterns 32-1, 31-2 and 21-3 and beginning with an increasing or decreasing pattern, by the complement operation. For instance, we have $N^{12\ldots k}_{13-2}(n) = N^{k(k-1)\ldots 1}_{31-2}(n)$.

\subsection{The pattern $12-3$} We first consider beginning with the pattern $p=k\dots 21$. In \cite[Lemma 9]{ClaesMans} it was proved that
   $$s_{12-3}(n;i)=\sum_{j=0}^{i-1}\binom{i-1}{j}s_{12-3}(n-2-j),$$
together with $s_{12-3}(n;n)=s_{12-3}(n;n-1)=s_{12-3}(n-1)$.

On the other hand, from the definitions, it is easy to see that
$$N_{12-3}^{k(k-1) \dots 1}(n)=\sum_{i=1}^{n-k+1}\binom{n-i}{k-1}s_{12-3}(n-k+1;i).$$
Hence, using (\ref{equality8}) and the fact shown in~\cite{Claes} that $s_{12-3}(n)$ equals $B_n$, we get the following proposition.

\begin{prop}
For all $n\geq k+1$, we have
$$\begin{array}{l}
N_{12-3}^{k(k-1) \dots 1}(n)=(k+1)B_{n-k}+\\
\qquad\qquad\qquad+\sum\limits_{j=0}^{n-k-2}
\left(\binom{n}{k+j}-k\binom{n-k-1}{j}-\binom{n-k}{j}\right)B_{n-k-1-j},
\end{array}$$ together with $N_{12-3}^{k(k-1) \dots 1}(k)=1$ and
$N_{12-3}^{k(k-1) \dots 1}(n)=0$ for all $n\leq k-1$.
\end{prop}

Now, let us consider beginning with the pattern $p=12\dots k$. From the definitions, it is easy to see that $N^{12\dots k}_{12-3}(n)=0$ for all $n$, where $k\geq
3$, and $N^{1}_{12-3}(n)=s_{12-3}(n)=B_n$ (see \cite[Prop. 10]{ClaesMans}). Thus, we only need to consider the case $k=2$.

Suppose $\pi\in S_{12-3}(n)$ is a permutation with $\pi_1<\pi_2$. It is easy to
see that $\pi_2=n$. Hence $N^{12}_{12-3}(n)=(n-1)s_{12-3}(n-2)$, for all
$n\geq 2$, and by \cite[Prop. 10]{ClaesMans}, we get the truth of the following

\begin{prop} We have

\[ E^{12\dots k}_{12-3}(x)=\left\{\begin{array}{ll}
0, & \mbox{if $k\geq3$,} \\[2mm]
x^2\displaystyle\sum_{j=0}^k(1-jx)^{-1}\displaystyle\sum_{d\geq0}\frac{x^d}{(1-x)(1-2x)\dots(1-dx)}, & \mbox{if $k=2$,} \\[4mm]
\displaystyle\sum_{d\geq 0}\frac{x^d}{(1-x)(1-2x)\dots(1-dx)}, & \mbox{if $k=1$.}
\end{array}\right.\]
\end{prop}

\subsection{The pattern $13-2$}
Let us introduce an object that plays an important role in the
proof of the main result in this case. For $n\geq m+1\geq 0$, we
define
$$A(n;m)=\displaystyle\sum_{1\leq i_m<\dots<i_2<i_1<n-1}s_{1-3-2}(n;i_1,i_2,\dots,i_m).$$
We extend this definition to $m=0$ by $A(n;0)=s_{1-3-2}(n)$.

\begin{lemma}\label{case132l}
For all $n\geq m\geq 0$,
$$A(n;m)=\displaystyle\sum_{j\geq0} (-1)^j\binom{m+1-j}{j}s_{1-3-2}(n-j).$$
\end{lemma}
\begin{proof}
For $m=0$ the lemma holds by definitions. Let $m\geq 0$; so
$$\begin{array}{ll}
A(n;m) &=\displaystyle\sum_{1\leq i_m<\dots<i_2<i_1<n-1}\sum_{j=1}^n s_{1-3-2}(n;i_1,i_2,\dots,i_m,j),\\
       &=A(n;m+1)+\displaystyle\sum_{1\leq i_m<\dots<i_2<i_1<n-1} s_{1-3-2}(n;i_1,i_2,\dots,i_m,n),\\
       &=A(n;m+1)+\displaystyle\sum_{1\leq i_m<\dots<i_2<i_1<n-1} s_{1-3-2}(n-1;i_1,i_2,\dots,i_m),\\
       &=A(n;m+1)+\displaystyle\sum_{1\leq i_m<\dots<i_2<n-2} s_{1-3-2}(n-1;n-1,i_2,\dots,i_m)+\\
       &\qquad\qquad\qquad\qquad+\displaystyle\sum_{1\leq i_m<\dots<i_2<i_1<n-2} s_{1-3-2}(n-1;i_1,i_2,\dots,i_m),\\
       &=A(n;m+1)+A(n-1;m)+\\
       &\qquad\qquad+\displaystyle\sum_{1\leq i_{m-1}<\dots<i_1<n-2} s_{1-3-2}(n-2;i_1,\dots,i_{m-1}),\\
       &=\cdots=A(n;m+1)+A(n-1;m)+\cdots+A(n-m-1;0).
\end{array}$$
Hence, using induction on $m$, we get

$$\begin{array}{l}
 A(n;m+1)=\displaystyle\sum_{j\geq0}(-1)^j\binom{m+1-j}{j}s_{1-3-2}(n-j)\\

\qquad\qquad\qquad\qquad-\displaystyle\sum_{d=0}^m\sum_{j\geq0}(-1)^j\binom{
m-d+1-j}{j}s_{1-3-2}(n-1-d-j).\end{array}$$

Using the identity
$\binom{r}{0}-\binom{r}{1}+\cdots+(-1)^s\binom{r}{s}=\binom{r-1}{s}$,
we get
$$\begin{array}{l}
 A(n;m+1)=\displaystyle\sum_{j\geq0}(-1)^j\binom{m+1-j}{j}s_{1-3-2}(n-j)\\
\qquad\qquad\qquad\qquad-\displaystyle\sum_{d=0}^m(-1)^d\binom{m-d}{d}s_{1-3
-2}(n-1-d).\end{array}$$
Now using the identity $\binom{n}{k}+\binom{n}{k-1}=\binom{n+1}{k}$, we get
$$A(n;m+1)=\displaystyle\sum_{j\geq0}(-1)^j\binom{m+2-j}{j}s_{1-3-2}(n-j),$$
which means that the lemma holds for $m+1$.
\end{proof}

Now we find $N^{k(k-1) \dots 1}_{13-2}(n)$.

\begin{prop}
Let $k\geq 1$. For all $n\geq 0$,
$$N^{k(k-1) \dots 1}_{13-2}(n)=C_{n+1-k}+\sum_{d=0}^{k-2}\sum_{j\geq0}(-1)^j\binom{k+1-d-j}{j}C_{n-d-j}.$$
\end{prop}

\begin{proof}
Claesson \cite{Claes} proved that the set of permutations that avoid
the pattern $13-2$ is the same as the set of permutations that avoid the pattern $1-3-2$, hence
\begin{equation}
N^{k(k-1) \dots 1}_{13-2}(n)=N^{k(k-1) \dots 1}_{1-3-2}(n).
\label{equality7}
\end{equation}

If the leftmost letter of a permutation avoiding 13-2 and beginning with the pattern $k(k-1)\ldots 1$ is $n$, then, obviously, there are $N^{(k-1)(k-2) \dots 1}_{1-3-2}(n-1)$ such permutations. Otherwise, it is easy to see that there are $A(n;k)$ such permutations. So, by Lemma~\ref{case132l} and the considerations above, also using the fact that the number of $(1-3-2)$-avoiding $n$-permutations in $S_n$ is $C_n$, we get
$$N^{k(k-1) \dots 1}_{13-2}(n)=N^{(k-1)(k-2) \dots 1}_{13-2}(n-1)+\displaystyle\sum_{j\geq0}(-1)^j\binom{k+1-j}{j}C_{n-j}.$$
Moreover, using the definitions and Equation~(\ref{equality7}), we have
$N^{1}_{13-2}(n)=s_{1-3-2}(n)=C_n$, for all $n\geq 0$. Hence, by induction on $k$, the proposition holds.
\end{proof}

Now, let us consider the case of $N^{12\dots k}_{13-2}(n)$.

\begin{prop}
Let $k\geq 1$. For all $n\geq k$, we have
$$N^{12\dots k}_{13-2}(n)=C_{n+1-k}.$$
\end{prop}

\begin{proof}
Suppose $\pi={\pi}^{\prime}n{\pi}^{\prime\prime}$ is a permutation in
$S_n(13-2)=S_n(1-3-2)$ (see~(\ref{equality7})), such that
$\pi_1<\pi_2<\dots<\pi_k$. It is easy to see that there exists an $m$ such that
    $$\pi=(m+1)(m+2)\dots(m+k-1)\beta n {\pi}^{\prime\prime},$$
where $\beta$ is a $1-3-2$-avoiding permutation on the letters $m+k,m+k+1,\dots,n-1$, and ${\pi}^{\prime\prime}\in S_{m}(1-3-2)$. Hence, in
terms of generating functions, we get
$$\sum_{n\geq0}N^{12\dots k}_{13-2}(n)x^n=x^kC^2(x).$$
The rest is easy to check using the identity $xC^2(x)=C(x)-1$.
\end{proof}

\subsection{The pattern $23-1$}
We first consider beginning with the pattern $p=k(k-1) \dots 1$.

\begin{prop}
For all $k\geq 1$,
$$E^{k(k-1) \dots 1}_{23-1}(x)=x^{k-1}\left(\sum_{d\geq0}\frac{x^d}{(1-x)(1-2x)\cdots(1-dx)}-1\right).$$
\end{prop}

\begin{proof}
Let $\pi\in S_n(23-1)$ be a permutation such that $\pi_1<\pi_2<\dots<\pi_k$. Since $\pi$ avoids $23-1$, we have $\pi_j=j$, for each $j=1,2,\dots,k-1$. Hence $\pi=12\dots (k-1){\pi}^{\prime}$, where $\pi^{\prime}$ is a non-empty $23-1$-avoiding permutation in $S_{n+1-k}$. The rest is easy to get by using \cite[Prop. 17]{ClaesMans}.
\end{proof}

Now let us consider beginning with the pattern $p=12\dots k$.

\begin{prop}
Suppose $k\geq 1$. For all $n\geq k+1$,
$$N^{12 \dots k}_{23-1}(n)=
\left(1+\binom{n-1}{k-1}\right)B_{n-k}+
\sum_{j=0}^{n-k-2}\left[\binom{n-1}{k+j}-\binom{n-k-1}{j}\right]B_{n-k-1-j},$$
with $N^{12 \dots k}_{23-1}(k)=1$.
\end{prop}
\begin{proof}
In \cite[Lemma 16]{ClaesMans} proved that for all $2\leq i\leq n-1$,
$$s_{23-1}(n;i)=\sum_{j=0}^{i-2}\binom{i-2}{j}s_{23-1}(n-2-j),$$
together with $s_{23-1}(n;n)=s_{23-1}(n;1)=s_{23-1}(n-1)=B_{n-1}$.

On the other hand, by the definitions, it is easy to see that
$$N_{23-1}^{12 \dots k}(n)=\sum_{i=1}^{n-k+1}\binom{n-i}{k-1}s_{23-1}(n-k+1;i).$$
Hence, using~(\ref{equality8}) and the fact that \cite{Claes} $s_{23-1}(n)$ is
given by $B_n$, we get the desired result.
\end{proof}


\section{Avoiding a pattern xy-z and beginning with the pattern $(k-1)(k-2)\ldots 1k$ or $23\ldots k1$}

In this section we consider avoidance of one of the patterns $12-3$, $13-2$, $23-1$, $21-3$, $31-2$ and $13-2$ and beginning with the pattern $(k-1)(k-2)\ldots 1k$. The case when a permutation begins with the pattern $23\ldots k1$ and avoids a pattern $xy-z$ can be obtained then by the complement operation.

\subsection{Avoiding $12-3$ and beginning with $(k-1)(k-2)\ldots 1k$}

\begin{prop} We have
$$N^{(k-1)(k-2) \ldots 1k}_{12-3}(n)=\binom{n-1}{k-1}B_{n-k}.$$
\end{prop}

\begin{proof}
Suppose $\pi={\pi}^{\prime}n{\pi}^{\prime\prime}$ avoids the pattern $12-3$ and begins with the pattern $(k-1)(k-2)\ldots 1k$. We have that ${\pi}^{\prime}$ must be decreasing, since otherwise we have an occurrence of the pattern $12-3$ involving the letter $n$, and ${\pi}^{\prime\prime}$ must avoid $12-3$. Also, since $\pi$ begins with $(k-1)\ldots21k$, the length of ${\pi}^{\prime}$ is $k-1$. Hence, by
\cite{Claes} (the number of permutations in $S_n(12-3)$ is given by $B_n$), we have
$$N^{(k-1)(k-2)\ldots 1k}_{12-3}(n)=\binom{n-1}{k-1}B_{n-k}.$$
\end{proof}

\subsection{Avoiding $13-2$ and beginning with $(k-1)(k-2)\ldots 1k$}

By \cite{Claes}, a permutation $\pi$ avoids the pattern $13-2$ if and only if $\pi$ avoids $1-3-2$.

Suppose $\pi={\pi}^{\prime}n{\pi}^{\prime\prime}$ is an $n$-permutation avoiding $1-3-2$ and beginning with $(k-1)(k-2)\ldots 1k$.  Obviously, ${\pi}^{\prime}$ and ${\pi}^{\prime\prime}$ avoid $1-3-2$ and each letter of ${\pi}^{\prime}$ is greater than any letter of ${\pi}^{\prime\prime}$, since otherwise we have an occurrence of the pattern $1-3-2$ involving the letter $n$. Also, ${\pi}^{\prime}$ begins with the pattern $(k-1)(k-2)\ldots 1k$ or ${\pi}^{\prime}=(k-1)(k-2)\ldots 1$.

By \cite{Knuth}, the generating function for the number of permutations that avoid $1-3-2$ is $C(x)$, hence, using the considerations above,
$$G^{(k-1)(k-2)\dots 1k}_{13-2}(x)=xG^{(k-1)(k-2)\dots 1k}_{13-2}(x)C(x)+x^{k}C(x).$$
Therefore, by~(\ref{catalan}), we get the following.

\begin{prop} We have
$$G^{(k-1)(k-2) \dots 1k}_{13-2}(x)=x^kC^2(x).$$
Hence
$$N^{(k-1)(k-2)\ldots 1k}_{13-2}(n) = \left\{ \begin{array}{ll} C_{n-(k-1)}, & \mbox{if $n\geq k$} \\ 0, & \mbox{else.}\end{array}
\right.$$
\end{prop}

\subsection{Avoiding $21-3$ and beginning with $(k-1)(k-2)\ldots 1k$}

If $k\geq 3$ then, by the definitions, we have $N^{(k-1)(k-2) \ldots 1k}_{21-3}(n)=0$.

If $k=1$ then, by the definitions and~\cite{Claes}, we have $N^{1}_{21-3}(n)=B_n$.

Suppose $k=2$ and $\pi={\pi}^{\prime}n{\pi}^{\prime\prime}$ is an $n$-permutation avoiding the pattern $21-3$ and beginning with the pattern $(k-1)(k-2)\ldots 1k = 12$. It is easy to see that ${\pi}^{\prime}$ must be increasing, and the length of ${\pi}^{\prime}$ is at least $1$. Thus, using the fact that the number of permutations in $S_n(21-3)$ is given by $B_n$ (see \cite{Claes}), we have
\begin{equation}
N^{(k-1)(k-2)\ldots 1k}_{21-3}(n)=\displaystyle\sum_{j=1}^{n-1}\binom{n-1}{j}B_{n-1-j}.
\label{equality16}
\end{equation}
Since $B_n=\displaystyle\sum_{j=0}^{n-1}\binom{n-1}{j}B_{n-1-j}$, equality~(\ref{equality16}) gives that
$$N^{(k-1)(k-2)\ldots 1k}_{21-3}(n)=B_n-B_{n-1}.$$
Thus we have proved the following.

\begin{prop}
$$N^{(k-1)(k-2)\ldots 1k}_{21-3}(n) = \left\{ \begin{array}{ll} 0, & \mbox{if $k\geq 3$} \\ B_{n}-B_{n-1}, & \mbox{if $k=2$,} \\ B_n, & \mbox{if $k=1$.} \end{array}
\right.$$
\end{prop}

\subsection{Avoiding $23-1$ and beginning with $(k-1)(k-2)\ldots 1k$}

\begin{prop} We have
$$N^{(k-1)\ldots 1k}_{23-1}(n)
=\left\{ \begin{array}{ll}
B_{n-k}+\displaystyle\sum_{t=2}^{n-k+2}\binom{t+k-3}{k-2}\displaystyle\sum_{j=0}^{t-2}\binom{t-2}{j}B_{n-k-1-j},
& \mbox{if $k\geq 3$} \\ B_{n-1}, & \mbox{if $k=2$,} \\ B_n,  &
\mbox{if $k=1$.} \end{array} \right.$$
\end{prop}

\begin{proof}
Suppose $k=2$. We are interested in the permutations $\pi\in S_n(23-1)$ that begin with the pattern $12$. It is easy to see that $\pi_1=1$, hence
$B_{12}^{23-1}(n)=B_{n-1}$ for all $n\geq2$.

Suppose $k\geq3$. We recall that $s_{23-1}(n;t)$ is the number of permutations in
$S_n(23-1)$ having $t$ as the first letter. By \cite{ClaesMans}, $s(n;1)=B_{n-1}$ and for $t\geq 2$, we have
$$s(n;t)=\sum_{j=0}^{t-2}\binom{t-2}{j}B_{n-2-j}.$$

On the other hand, if a permutation $\pi={\pi}^{\prime}1t{\pi}^{\prime\prime}$ avoids $23-1$ and begins with the pattern $(k-1)(k-2)\ldots 1k$, then ${\pi}^{\prime}$ is decreasing of length $k-2$, and using $s(n;t)$, we get
$$N^{(k-1)(k-2)\ldots 1k}_{23-1}(n)=
B_{n-k}+\displaystyle\sum_{t=2}^{n-k+2}\binom{t+k-3}{k-2}\displaystyle\sum_{j=0}^{t-2}\binom{t-2}{j}B_{n-k-1-j}.$$
\end{proof}

\subsection{Avoiding $31-2$ and beginning with $(k-1)(k-2)\ldots 1k$}

By \cite{Claes}, a permutation $\pi$ avoids the pattern $31-2$ if and only if $\pi$ avoids the pattern $3-1-2$.

Suppose $\pi={\pi}^{\prime}1{\pi}^{\prime\prime}$ is an $n$-permutation avoiding $3-1-2$ and beginning with $(k-1)(k-2)\ldots 1k$.  Obviously, ${\pi}^{\prime}$ and ${\pi}^{\prime\prime}$ avoid $3-1-2$ and each letter of ${\pi}^{\prime}$ is smaller than any letter of ${\pi}^{\prime\prime}$, since otherwise we have an occurrence of the pattern $3-1-2$ involving the letter 1. Also, ${\pi}^{\prime}$ begins with the pattern $(k-1)(k-2)\ldots 1k$ or ${\pi}^{\prime}=(k-1)(k-2)\ldots 2$ and ${\pi}^{\prime\prime}$ is not empty. So, using the generating function for the number of
permutations avoiding the pattern $3-1-2$, which is $C(x)$ (\cite{Knuth}), we get
$$G^{(k-1)(k-2)\dots 1k}_{31-2}(x)=xG^{(k-1)(k-2)\dots 1k}_{31-2}(x)C(x)+x^{k-1}(C(x)-1).$$
Therefore, using~(\ref{catalan}), we get the following.

\begin{prop} We have
$$G^{(k-1)(k-2) \dots 1k}_{31-2}(x)= \left\{ \begin{array}{ll} x^kC^3(x), & \mbox{if $k \geq 2$,} \\ C(x), & \mbox{if $k=1$.} \end{array}
\right.$$
Hence
$$N^{(k-1)(k-2)\ldots 1k}_{31-2}(n)= \left\{ \begin{array}{ll} C_{n-k+2}-C_{n-k+1}, & \mbox{if $k\geq 2$,} \\ C_n, & \mbox{if $k=1$.} \end{array}
\right.$$
\end{prop}

\subsection{Avoiding $32-1$ and beginning with $(k-1)(k-2)\ldots 1k$}

\begin{prop}
$$N^{(k-1)(k-2)\ldots 1k}_{32-1}(n) = \left\{ \begin{array}{ll} 0, & \mbox{if $k\geq 4$} \\ B_{n-1}-(n-2)B_{n-3}, & \mbox{if $k=3$ and $n\geq3$,} \\ B_{n}-(n-1)B_{n-2}, & \mbox{if $k=2$ and $n\geq 2$,} \\ B_n,  & \mbox{if $k=1$.} \end{array}
\right.$$
\end{prop}

\begin{proof}
Using the definitions and \cite{Claes}, it is easy to see that the statement is true for $k=1,2$ and $k\geq 4$.

Suppose now that $k=3$ and $\pi={\pi}^{\prime}1{\pi}^{\prime\prime}$ is an $n$-permutation avoiding the pattern $32-1$ and beginning with the pattern $(k-1)(k-2)\ldots 1k = 213$. We have that ${\pi}^{\prime}$ must be increasing, since otherwise we have an occurrence of the pattern $32-1$ involving the letter $1$, and ${\pi}^{\prime\prime}$ must avoid $32-1$. Moreover, since $\pi$ begins with $213$, the length of $\pi$ is $1$ and the rightmost letter of ${\pi}^{\prime\prime}$ is greater than the letter of ${\pi}^{\prime}$. Also, it is easy to see that the
number of permutations in $S_{n-1}(32-1)$ beginning with the pattern $12$ is
the same as the number of permutations in $S_n(32-1)$ beginning with the pattern $213$
(one can see it by placing $1$ in the second position). Hence $N^{(k-1)\ldots21k}_{32-1}(n)=B_{n-1}-(n-2)B_{n-3}$ for all $n\geq3$.
\end{proof}

\section{Avoiding a pattern x-yz and beginning with the pattern $(k-1)(k-2)\ldots 1k$ or $23\ldots k1$}

In this section we consider avoidance of one of the patterns $1-23$, $1-32$, $2-31$, $2-13$, $3-12$ and $1-32$ and beginning with the pattern $(k-1)(k-2)\ldots 1k$. The case when a permutation begins with the pattern $23\ldots k1$ and avoids a pattern $x-yz$ can be obtained by the complement operation.

\begin{prop} We have
\[ E_{1-32}^{(k-1)(k-2)\ldots 1k}(x) = \left\{ \begin{array}{ll} e^{e^x}\int_{0}^{x}e^{-e^t}\displaystyle\sum_{n\geq k-1}\frac{t^n}{n!}\ dt, & \mbox{if $k\geq 2$,}
\\[2mm] e^{e^x-1}, & \mbox{if $k=1$.}
\end{array}
\right. \]
\end{prop}

\begin{proof} Suppose $k\geq 2$. Let $B_{n, k}$ denote the number of $n$--permutations that
avoid the pattern $1-32$ and begin with the pattern $(k-1)(k-2)\ldots 1k$. Suppose $\pi= \sigma 1 \tau$ is such a permutation of length $n+1$. Obviously, the letters of $\tau$ must be in increasing order, since otherwise we have an occurrence of the pattern $1-32$ in $\pi$ starting from the letter $1$. If $|\sigma| = i$, then we can choose the letters of $\sigma$ in ${n \choose i}$ ways. Since the letters of $\tau$ are in increasing order, they do not affect $\sigma$ and thus there are $B_{i, k}$ possibilities to choose $\sigma$. Also, if $n\geq k-1$, then $1$ can be in the $(k-1)$th position, and in this case, since $\pi$ begins with the pattern $(k-1)(k-2)\ldots 1k$, it must be that $\pi = (k-1)(k-2)\ldots 21k(k+1)\ldots (n+1)$. Thus, in the last case we have only one permutation. This leads to the recurrence relation
$$B_{n+1,k} = \displaystyle\sum_{i\geq 0}{n \choose i}B_{i,k} + \delta_{n,k},$$
where
\[ \delta_{n,k} = \left\{ \begin{array}{ll} 1, & \mbox{if $n\geq k-1$,}
\\ 0, & \mbox{else.}
\end{array}
\right. \]

This recurrence relation is identical to the one given in the proof of Proposition~\ref{about_3-12}, so using this proof we get the desired result. \end{proof}

\begin{prop} We have
\[ E_{1-23}^{(k-1)(k-2)\ldots 1k}(x) = \left\{ \begin{array}{ll} e^{e^x}\int_{0}^{x}\int_{0}^{t}\frac{r^{k-2}}{(k-2)!}e^{r-e^t}\ drdt, & \mbox{if $k\geq 2$,}
\\[2mm] e^{e^x-1}, & \mbox{if $k=1$.}
\end{array}
\right. \]
\end{prop}

\begin{proof} If $k=1$, then the statement is true due to Proposition~\ref{about_1-23}.

Suppose $k\geq 2$. Let $B_{n, k}$ denote the number of $n$-permutations that
avoid the pattern $1-23$ and begin with the pattern $(k-1)(k-2)\ldots 1k$. Suppose $\pi= \sigma 1 \tau$ is such a permutation of length $n+1$. Obviously, the letters of $\tau$ must be in decreasing order since otherwise we have an occurrence of the pattern $1-23$ in $\pi$ starting from the letter $1$. If $|\sigma| = i$, then we can choose the letters of $\sigma$ in ${n \choose i}$ ways. Since the letters of $\tau$ are in the decreasing order, they do not affect $\sigma$ and thus there are $B_{i, k}$ possibilities to choose $\sigma$. Besides, if $n\geq k-1$, then $1$ can be in the $(k-1)$th position, and in this case, since $\pi$ begins with the pattern $(k-1)(k-2)\ldots 1k$ and $\tau$ is decreasing, it must be that the $k$th letter of $\pi$ is $(n+1)$ and there are ${n-1 \choose k-2}$ ways to choose the letters of $\sigma$ and then write them in decreasing order. Thus,
$$B_{n+1,k} = \displaystyle\sum_{i\geq 0}{n \choose i}B_{i,k} + {n-1 \choose k-2}.$$

Multiplying both sides of the equality with $x^n/n!$ and summing over $n$, we get the differential equation
$$\frac{d}{dx}E_{1-23}^{(k-1)(k-2)\ldots 1k}(x) = E_{1-23}^{(k-1)(k-2)\ldots 1k}e^x + \displaystyle\sum_{n\geq 0}{n-1 \choose k-2}\frac{x^n}{n!},$$
with the initial condition $E_{1-23}^{(k-1)(k-2)\ldots 1k}(0)=0$. If $F(x)$ denotes the last term, then it is easy to see that $F^{\prime}(x) = \frac{x^{k-2}}{(k-2)!}e^x$, and thus
$$F(x)=\int_{0}^{x}\frac{t^{k-2}}{(k-2)!}e^t\ dt.$$

Now, the solution to the equation above is given by
\begin{equation}
E_{1-23}^{(k-1)(k-2)\ldots 1k}(x)= e^{e^x}\int_{0}^{x}e^{-e^t}F(t)\ dt = e^{e^x}\int_{0}^{x}\int_{0}^{t}\frac{r^{k-2}}{(k-2)!}e^{r-e^t}\ drdt.
\label{equality18}
\end{equation}

For example, if $k=2$, then $(k-1)(k-2)\ldots 1k = 12$ and (\ref{equality18}) gives
$$E_{1-23}^{12} = e^{e^x}\int_{0}^{x}e^{-e^t}(e^t - 1)\ dt,$$
which is a particular case of Proposition~\ref{about_3-21}, since the number of $n$-permutations that avoid the pattern $3-21$ and begin with the pattern $21$ is equal to the number of $n$-permutations that avoid the pattern $1-23$ and begin with the pattern $12$ by applying the complement.
\end{proof}

\begin{prop} We have
$$G^{(k-1)(k-2) \ldots 1k}_{2-13}(x)= \left\{ \begin{array}{ll} 0, & \mbox{if $k\geq 3$} \\ x^2C^3(x), & \mbox{if $k=2$} \\ C(x), & \mbox{if $k=1$.}\end{array}
\right.$$
Hence
$$N^{(k-1)(k-2)\ldots 1k}_{2-13}(n)=\left\{ \begin{array}{ll} 0, & \mbox{if $k\geq 3$} \\ C_{n-1}-C_{n-2}, & \mbox{if $k=2$} \\ C_n, & \mbox{if $k=1$.}\end{array}
\right.$$
\end{prop}
\begin{proof}
For the case $k=1$, see Proposition~\ref{about_2-13}. If $k\geq 3$, then the statement is true, since in this case the pattern $(k-1)(k-2) \ldots 1k$ does not avoid $2-13$.

Suppose now that $k=2$. Using the reverse, we see that beginning with the pattern $12$ and avoiding $2-13$ is equivalent to ending with the pattern $21$ and avoiding $31-2$, which by \cite{Claes} is equivalent to ending with the pattern $21$ and avoiding the pattern $3-1-2$.

Let $\pi={\pi}^{\prime}1{\pi}^{\prime\prime}$ be an $n$-permutation avoiding $3-1-2$ and ending with the pattern $21$.  Obviously, ${\pi}^{\prime}$ and ${\pi}^{\prime\prime}$ avoid $3-1-2$ and each letter of ${\pi}^{\prime}$ is less than any letter of ${\pi}^{\prime\prime}$, since otherwise we have an occurrence of $3-1-2$ involving the letter 1. Also, ${\pi}^{\prime\prime}$ ends with the pattern $21$ or $|{\pi}^{\prime\prime}| = 1$. So, using the generating function for the number of permutations avoiding $3-1-2$, which is $C(x)$ (\cite{Knuth}), we have
$$G^{12}_{2-13}(x)=xG^{12}_{2-13}(x)C(x)+x(C(x)-1).$$

Therefore, using~(\ref{catalan}), we get the desired result. \end{proof}

\begin{prop} We have
$$G^{(k-1)(k-2) \dots 1k}_{2-31}(x)=x^kC^2(x).$$
Hence
$$N^{(k-1)(k-2)\ldots 1k}_{2-31}(n) = \left\{ \begin{array}{ll} C_{n-(k-1)}, & \mbox{if $n\geq k$} \\ 0, & \mbox{else.}\end{array}
\right.$$
\end{prop}

\begin{proof}
Using the reverse, we see that beginning with the pattern $(k-1)(k-2)\ldots 1k$ and avoiding the pattern $2-31$ is equivalent to ending with the pattern $k12\ldots (k-1)$ and avoiding the pattern $13-2$, which, by \cite{Claes}, is equivalent to ending with the pattern $k12\ldots (k-1)$ and avoiding the pattern $1-3-2$.

Let $\pi={\pi}^{\prime}n{\pi}^{\prime\prime}$ be an $n$-permutation avoiding the pattern $1-3-2$ and ending with the pattern $k12\ldots (k-1)$.  Obviously, ${\pi}^{\prime}$ and ${\pi}^{\prime\prime}$ avoid the pattern $1-3-2$ and each letter of ${\pi}^{\prime}$ is greater than any letter of ${\pi}^{\prime\prime}$, since otherwise we have an occurrence of the pattern $1-3-2$ involving the letter $n$. Also, ${\pi}^{\prime\prime}$ ends with the pattern $k12\ldots (k-1)$ or ${\pi}^{\prime\prime}=12\ldots (k-1)$.

Using the reverse operation, the generating function for the number of permutations ending with the pattern $k12\ldots (k-1)$ and avoiding $1-3-2$ is equal to $G^{(k-1)(k-2)\ldots 1k}_{2-31}(x)$. In terms of generating functions, the considerations above lead to
$$G^{(k-1)(k-2)\dots 1k}_{2-31}(x)=xG^{(k-1)(k-2)\dots 1k}_{2-31}(x)C(x)+x^{k}C(x).$$

Therefore, by~(\ref{catalan}), we get the desired result.
\end{proof}

\begin{prop}\label{aaa} We have
\[ E_{3-12}^{(k-1)(k-2)\ldots 1k}(x) = \left\{ \begin{array}{ll} (e^{e^x}/(k-1)!)\int_{0}^{x}t^{k-1}e^{-e^t+t}\ dt, & \mbox{if $k\geq 2$,}
\\[2mm] e^{e^x-1}, & \mbox{if $k=1$.}
\end{array}
\right. \]
\end{prop}

\begin{proof} Suppose $k\geq 2$. Let $B_{n, k}$ denote the number of $n$-permutations that
avoid the pattern $3-12$ and begin with a decreasing subword of
length $k$. Let $\pi = \sigma (n+1) \tau$ be such a permutation of length $n+1$. Obviously, the letters of $\tau$ must be in decreasing order since otherwise we have an occurrence of $3-12$ in $\pi$ starting from the letter (n+1). If $|\sigma| = i$ then we can choose the letters of $\sigma$ in ${n \choose i}$ ways. Since the letters of $\tau$ are in decreasing order, they do not affect $\sigma$ and thus there are $B_{i, k}$ possibilities to choose $\sigma$. Also, if $|\sigma| = k-1$ and the letters of $\sigma$ are in decreasing order, we get ${n \choose k-1}$ additional ways to choose $\pi$. Thus
$$B_{n+1,k} = \displaystyle\sum_{i\geq 0}{n \choose i}B_{i,k} + {n \choose k-1}.$$
This recurrence relation is identical to the one given in the proof of Proposition~\ref{about_1-23}, and we get the desired result using that proof.
\end{proof}

\begin{prop}
$$E^{(k-1)(k-2)\ldots 1k}_{3-21}(n) = \left\{ \begin{array}{ll} 0, & \mbox{if $k\geq 4$} \\ (e^{e^x}/(k-1)!)\int_{0}^{x}t^{k-1}e^{-e^t+t}\ dt, & \mbox{if $k=2$ or $k=3$,} \\ e^{e^x-1}, & \mbox{if $k=1$.} \end{array}
\right.$$
\end{prop}

\begin{proof}
If $k\geq 4$ then the statement is true, since in this case the pattern $(k-1)(k-2)\ldots 1k$ does not avoid the pattern $3-21$. In the other cases, we use the same arguments as we have in the proof of Proposition~\ref{aaa}. The only difference is that instead of decreasing order in $\tau$, we have increasing order.
\end{proof}

\section{Conclusions}

The goal of our paper is to give a complete description for the numbers of permutations avoiding a pattern of the form $x-yz$ or $xy-z$ and either beginning with one of the patterns $12\ldots k$, $k(k-1)\ldots 1$, $23\ldots k1$, $(k-1)(k-2)\ldots 1k$, or ending with one of the patterns $12\ldots k$, $k(k-1)\ldots 1$, $1k(k-1)\ldots 2$, $k12\ldots (k-1)$. This description is given in Sections~5--8. However, some of our results can be generalized to beginning with a pattern belonging to $\Gamma_k^{min}$ or $\Gamma_k^{max}$, and thus to the ending with a pattern belonging to $\Delta_k^{min}$ or $\Delta_k^{max}$ (see Section~\ref{Preliminaries} for definitions). An example of such a generalisation is given in Theorem~\ref{generalization_2} below. This theorem generalizes Propositions~\ref{about_1-23} and \ref{aaa} and can be proved by using the same considerations as we do in the proofs of these propositions.

\begin{thm}\label{generalization_2} Suppose $p_1, p_2 \in \Gamma_k^{min}$ and $p_1 \in S_k(1-23)$, $p_2 \in S_k(1-32)$. Thus, the complements $C(p_1), C(p_2) \in \Gamma_k^{max}$ and $C(p_1) \in S_k(1-23)$, $C(p_2) \in S_k(3-12)$. Then, we have
$$E^{p_1}_{1-23}(x) = E^{C(p_1)}_{3-21}(x) = E^{p_2}_{1-32}(x) = E^{C(p_2)}_{3-12}(x) =$$
$$\left\{ \begin{array}{ll} (e^{e^x}/(k-1)!)\int_{0}^{x}t^{k-1}e^{-e^t+t}\ dt, & \mbox{if $k\geq 2$,}
\\ e^{e^x-1}, & \mbox{if $k=1$.}
\end{array}
\right.$$
\end{thm}

{\bf Acknowledgments}: The authors are grateful to Einar
Steing\'{\i}msson for his helpful suggestions and comments.

\end{document}